\documentclass[12pt,a4paper,twoside,reqno]{amsart}
\usepackage{amsmath}
\usepackage{amsfonts}
\usepackage{amssymb}
\usepackage{amstext}
\usepackage{amsbsy}

\usepackage{epsf}
\newcommand{\Proof}{\noindent \textbf{Proof} \par}

\theoremstyle{plain}
\begingroup

\newtheorem{thm}{Theorem}
\newtheorem{prop}[thm]{Proposition}
\newtheorem{cor}[thm]{\textbf{\textit{Corollary}}}
\newtheorem{lem}[thm]{Lemma}

\endgroup
\newtheorem{rem}[thm]{\textbf{Remark}}

\theoremstyle{definition}
\newtheorem{defn}[thm]{Definition}

\newtheorem{exa}[thm]{Example}

\newcommand{\ch}{\cosh}
\newcommand{\sh}{\sinh}
\newcommand{\ta}{\tanh}

\newcommand{\argch}{\text{arcosh}}
\def\o{\omega}

\newcommand{\ov}[1]{\overline{#1}}
\newcommand{\wt}[1]{\widetilde{#1}}

\newcommand{\rme}{\text{\rm Re}}
\newcommand{\rmi}{\text{\rm Im}}
\newcommand{\Om}{\Omega}
\newcommand{\hi}[1]{\mathbb{H}^#1}
\newcommand{\m}[1]{\mathbb{R}^#1}
\newcommand{\co}{\mathbb{C}}

\newcommand{\Sp}{\mathbb{S}}
\newcommand{\M}{\mathbb{M}}
\newcommand{\zed}{\mathbb{Z}}
\newcommand{\R}{\mathbb{R}}
\newcommand{\sd}{\mathbb{S}^2}

\newcommand{\hd}{\mathbb{H}^2}

\newcommand{\C}{\mathbb{C}}

\newcommand{\di}{\mathbb{D}}

\newcommand{\mk}{\R^{2,1}}
\newcommand{\si}{\sigma}
\newcommand{\f}{h}

\newcommand{\te}{\theta}
\newcommand{\ga}{\gamma}
\newcommand{\al}{\alpha}
\newcommand{\be}{\beta}

\def\rmd{\mathop{\rm d\kern -1pt}\nolimits}

\begin{document}
\centerline{\large\bf Associate and conjugate minimal immersions in $\M\times \R$}

\vskip4mm
\centerline{\textbf{\textit{Laurent Hauswirth}} {}
\textbf{\textit{Ricardo Sa Earp}} {}
\textbf{\textit{Eric Toubiana}}}
\title{ }

\address{Universit\'e de Marne-la-Vall\'ee, Marne-la-Vall\'ee, France\\
\newline email: Laurent.Hauswirth@univ-mlv.fr}

\address{Pontif\'\i cia Universidade Cat\'olica do Rio de janeiro\\
         Depto de matem\'atica\\
         Rua Marqu\^es de S\~ao Vicente 225\\
         24 453-900 Rio de Janeiro-RJ, Brazil\\
         fax:55-021-3114 1282\\
\newline email: earp@mat.puc-rio.br}

\address{Centre de Math\'ematiques de Jussieu, Universit\'e Paris VII\\
         Denis Diderot\\
         2 Place Jussieu\\
         F-75251 Paris Cedex 05, France\\
\newline email: toubiana@math.jussieu.fr}

\thanks{
The  authors likes to thanks CNPq and PRONEX of Brazil and Accord Brasil-France,
 for partial financial support.}

\date{\today}

\maketitle

\vskip8mm

\section{Introduction}\label{Sec.Intro}

A beautiful phenomenon in Euclidean space is the existence of a
$1$-parameter  family of minimal isometric surfaces connecting the
catenoid and the helicoid. They are {\em associate}.  A well-known
fact, is that any two conformal isometric minimal surfaces in
a space form are associate. What happens in other
$3$-dimensional manifolds ?

In this paper we will discuss the same phenomenon in the product
space, $\M\times \R,$  giving a definition of
associate minimal immersions.  We specialize in the situations $\M=\hd,$
the hyperbolic plane, and $\M=\mathbb S^2,$ the sphere where
surprising facts occur. We will
prove some existence and uniqueness results explained in the
sequel. We begin with the definition.

 Let $\M$ be a two dimensional Riemannian manifold. Let $(x, y,
t)$ be local  coordinates in $\M\times \R,$ where $z = x +iy$
are conformal coordinates on $\M$ and $t\in \R.$ Let $\sigma^2
|\rmd z|^2,$ be the conformal metric in $\M,$ hence $\rmd s^2 =
\sigma^2 |\rmd z|^2 +\rmd t^2$ is the metric in the product space
$\M\times \R.$ Let  $\Omega\subset \co$ be a planar simply connected
domain, $w =u+iv\in \Omega.$  We recall that if
$X:\Omega\rightarrow \M\times \R$,
$ w\mapsto (h(w), f (w)),$\newline
$ w\in \Omega$,
is a conformal minimal immersion
then
 $h:\Omega \rightarrow (\M, \sigma^2 |\rmd z|^2)$ is
a harmonic map.
 We recall also that for any harmonic map
$h :\Omega\subset\co \rightarrow \M$ there exists a related
  Hopf holomorphic function $Q (h)$. Two conformal immersions
$X=(h, f ),\  X^*=( h^* , f^* ) :\Omega \rightarrow \hd \times \R$
are  said {\em associate} if they are isometric
and if the Hopf
functions satisfy the relation $Q ( h^* ) =e^{2i\theta} Q(h)$
for a real number $\theta$. If
$Q ( h^* )= -Q (h)$ then the two immersions are said
 {\em conjugate}.

In this paper we will show that there exist two conformal
isometric minimal surfaces in $\hd\times \R,$ with constant
{\em Gaussian curvature} $-1$, that are non associate. We will
prove also that the vertical cylinder over a planar geodesic in
$\hd\times \R,$ are the only minimal surfaces with constant
{\em Gaussian curvature} $K\equiv 0.$

One of our principal results is a uniqueness theorem in $\hd\times
\R,$ or $\mathbb S^2\times \R,$ showing that the conformal
metric and the Hopf function determine a minimal conformal
immersion, up to an isometry of ambient space, see
Theorem \ref{T.immiso}. We will derive the
existence of the minimal associate family in $\hd\times \R$ and $\sd \times \R$ in corollary \ref{C.associate} by establishing an
existence result, see Theorem \ref{T.existence}.
The associate minimal family in $\hd\times \R$ and
$\sd \times \R$ is
derived by another approach in \cite{[Dan]}.

The first author has constructed examples of minimal surfaces in
$\hd \times \R$ and in $\sd \times \R$ which generalize the
family of Riemann's minimal examples of $\R^3$. He classify and construct all
example foliated by horizontal constant curvature curves. Some of them
have Gaussian curvature $K\equiv -1$. This family is parametrized
by two parameter  $(c,d)$ and the example corresponding to
$(c,d)$ is conjugate to the one parametrized by $(d,c)$ (we refer to the
paper \cite{[H]}  for more details on these surfaces).
The second and third authors, proved that any two minimal
isometric  screw motion immersions in $\hd\times \R$ and
$\mathbb S^2\times \R$ are associate, see \cite{[S-T]}.
The second author proved  that any
two minimal isometric parabolic screw motion immersions into
$\hd\times \R$ are associate. On the other hand,
he proved that there exist families of
{\em associate hyperbolic} screw motion immersions, but there exist also
isometric non-associate {\em hyperbolic} screw motion immersions,
see \cite{[Earp]}. There exists hyperbolic screw motion
surfaces associate to parabolic screw motion surfaces
(see example 13).

Several questions arise from this work: We point out the
problem of the existence of the associate minimal family in
$\M\times \R,$ for any $2$-dimensional Riemannian manifold $\M$.
Also we may ask in which general assumptions isometric immersions must
be associate ?

The second principal result is a generalization of the Krust's theorem (see \cite{hil}, tome I,  page 118 
and applications therein) which state that an associate surface
of a minimal vertical graph on a convex domain is a vertical graph. These theorem is true in
$\M \times \R$ when the Gaussian curvature $K_{\M} \leq 0$. It will not be
true anymore by example in  $\sd \times \R$ (see Theorem \ref{krust}).

For related works on minimal surfaces in $\M\times \R$, see for
instance Daniel \cite{[Dan]} Nelli and Rosenberg \cite{[N-R]}, Meeks and Rosenberg
\cite{[M-R]} and Rosenberg \cite{[R]}.

\bigskip

\section{Preliminar} \label{Sec.Prel}
We consider $X: \Om \subset \mathbb{R} ^2 \rightarrow \M \times \R$
a minimal surface conformaly embedded in a product space. $\M$ is a Riemannian complete
two-manifold with metric $\mu=\si^2 (z) \vert dz \vert ^2 $ and Gauss curvature $K_\M$. First we fix some notations. Let us denote $|v|^2_\si=\si^2 | v|$, 
$\langle v_{1},v_{2} \rangle_{\si}=\si ^2 \langle v_{1},v_{2} \rangle$ where $|v|$ and $<v_{1},v_{2}>$ stands for the standard norm and inner product in $\R ^2$.
Let us find $w=u+iv$ as  conformal parameters of $\Om$, i.e.
$ds_{X}^2= \lambda ^2 \vert \rmd w \vert^2$. We denote by $X=(\f,f)$ the immersion
where $\f(w) \in \M$ and $f(w) \in \R$.

Assume that $\M$ is isometrically embedded in $\R^k$. 
By definition (see Lawson \cite{lawson}) the mean curvature vector in is
$$2\overrightarrow{H}=(\triangle X ) ^{T_X \M \times \R}=((\triangle \f)^{T_\f
\M},\triangle f)=0$$
\noindent
where $h=(h_{1},...,h_{k})$. Then $\f : \Om \longrightarrow \M$ is a harmonic map
between $\Om$ and the complete Riemannian surface $\M$ and $f$ is a real harmonic function. 
If $(U,  \si^2 (z) \vert \rmd z \vert^2)$ is a local
parametrization of $\M$, the harmonic map  equation
in the complex coordinate $z=x+i y$ of $\M$ (see \cite{sy}, page 8) is

\begin{equation}
\f_{w \bar{w}}+ 2(\log \si \circ \f )_z \f_w \f_{\bar{w}}=0
\label{eq:1.1}
\end{equation}

In the theory of {\it harmonic map} there is two important
classical object to consider. One is the holomorphic {\it quadratic
Hopf differential} associate to $\f$:
\begin{equation}
Q(\f)=(\si \circ h)^2 h_w \ov{h}_w (\rmd w)^{2}:= \phi (w) (dw)^2
\label{eq:hopf}
\end{equation}

\noindent
An other object is the {\it complex coefficient of dilatation} (see Alhfors \cite{alhfors})
of a quasi-conformal map:

\begin{equation*}
a(w)=\frac{\overline{\f_{\bar w}}}{\f_{w}}
\end{equation*}

Since we consider conformal immersion, we have
$(f_w)^{2}=-\phi (w)$ from (see \cite{[S-T]}):
\[
\begin{array}{l}
\vert \f _u \vert^2_{\si} + (f_u)^2  = \vert \f _v \vert ^2_{\si} + (f_v)^2 \\
\\
\left<\f_u ,\f _v \right>_{\si}+ f_u.f_v =0
\end{array}
\]

\noindent
We define $\eta$ as the holomorphic
one form $\eta =\pm 2i \sqrt{\phi (w)}\rmd w$ when $\phi$ have
only even zeroes. The sign is chosen in function of $f$ to
have:

\begin{equation}
f= {\rm Re} \int _w \eta
\label{eq:third}
\end{equation}

 In the case where $\M=\R^2$, $g$ is the Gauss map in the
classical Weierstrass representation.

When $X$ is a conformal immersion then the Gauss map $N$ in
$\M ^2 \times \R$ is
 given by (see \cite{[S-T]}):
\begin{equation}\label{Eq.Gaussmap}
N=\dfrac{(\frac{2}{\sigma}\rme g, \frac{2}
{\sigma}\rmi g,\vert  g\vert^2-1)}{\vert g\vert^2+1}
\end{equation}
where  $g :=\dfrac{ f _wh_{\ov{w}}- f _{\ov{w}}h_w} {\sigma\vert h_{\ov{w}}\vert (\vert h_w \vert + \vert h_{\ov{w}}\vert)}$.
We remark that

 \begin{equation}
 g^{2}=-\dfrac{h_w}{\ov{h _{\bar w}}}
\label{eq:dilatation2}
\end{equation}

Using the equations above (\ref{eq:hopf}),(\ref{eq:dilatation2}) we
can express the differential $dh$ as follows:

\begin{equation}
d \f= \f_{\bar w} d\bar w+\f_{w} dw=
\frac{1}{2 \si}\overline{ g^{-1} \eta} -\frac{1}{2
\si} g \eta
\label{eq:1}
\end{equation}

\noindent
The metric $d s^2_X=\lambda ^2 \vert \rmd w\vert ^2$ is given
by \cite{[S-T]}:

\begin{equation}
\label{eq:metricbis} ds^2_X=(\vert h_{w} \vert_\si + \vert h_{\bar
w} \vert_\si ) ^2 \vert \rmd w\vert ^2
\end{equation}

Thus combining together the equation we derive the metric in terms
of $g$ and $\eta$:

\begin{equation}
ds^{2}_{X}=\frac14(\vert g\vert^{-1}+\vert g \vert)^{2}\vert \eta \vert ^{2}
=(\vert\sqrt{a}\vert+\vert \sqrt{a} \vert^{-1})^{2}\vert \phi \vert \vert dw \vert^{2}
\label{eq:metric}
\end{equation}

In the case of minimal surfaces $X$ conformally immersed in
$\R^3=\R^2 \times \R$, the data $(g,\eta)$ are classical
Weierstrass data:

$$X(w)=(\f, f)=\left(\frac12\overline{\int_{w} g^{-1} \eta} -\frac12
\int_{w} g \eta ,{\rm Re} \int_{w} \eta \right)$$

The main difference is in the fact that $g$ is no more a meromorphic map in the case where the ambient space
is $\M \times \R$.
To study $g$ it is more
convenient for some purpose
to consider the complex function $\o + i \psi$
defined by

\begin{equation}
\label{eq:gauss}
g:=-i e^{\o+i\psi}
\end{equation}

We will derive some equation for $\psi$ in Lemma \ref{cauchy} in
term of $\o$. It is a well known fact (see \cite{sy} page 9) that
harmonic mappings satisfy the B\"ochner formula:

\begin{equation}
\triangle _0 \log \frac{ \vert \f_w   \vert} {\vert \f _
{\bar w} \vert } = -2K_\M J(\f)
\label{eq:bochner}
\end{equation}

\noindent where $\displaystyle {J(\f)= \si^2 \left(  \vert \f_w
\vert ^2 -\vert \f _ {\bar w} \vert^2 \right)}$ is the Jacobian of
$\f$ with $\vert \f_w  \vert^2= \f _w \overline{\f_w}$. Hence
taking into account (\ref{eq:hopf}), (\ref{eq:dilatation2}), (\ref{eq:gauss}) and
(\ref{eq:bochner}):
\begin{equation}
\triangle _0 \o = -2K_\M \sh(2\o)\vert \phi \vert
\label{eq:gordon}
\end{equation}
\noindent where $\triangle_{0}$ denote the laplacian in the
euclidean metric. With these convention notice that the metric and
the third coordinate of the Gauss map $N$ are given by

$$ds^2_X=4\cosh ^2 \o \vert \phi \vert \vert \rmd w\vert ^2
\;\; \hbox{ and } \;\; N_3 = \ta \o $$

\bigskip
On account of the above discussion we deduce the following:

\begin{prop}\label{P.min.imm}
Let $h:\Omega\rightarrow \M$ be a harmonic mapping such that
the holomorphic quadratic differential $Q (h)$ does not vanish or so have zero with
even order. Then there exists a complex map $g=-ie^{\o+i\psi}$
and a holomorphic one form $\eta = \pm 2i \sqrt Q(h)$ such that, with $f= \rme \int \eta$,
the map $X:=(h, f ): \Omega\rightarrow \M\times \R$ is a conformal and minimal (possibly branched) immersion.
The third component of the normal vector is given by $N_3=\ta \o$. The metric of the immersion is given by (\ref{eq:metric}):

$$ds^2_X= \cosh ^2 \o \vert \eta \vert^2$$
\noindent where $\o$ is a solution of the sh-Gordon equation

$$\triangle_0 \o =- 2K_\M \sinh(2\o)\vert \phi \vert.$$

\end{prop}

\Proof
We deduce from the hypothesis that we can
solve in $ f $ the equation
$( f _w)^2=-(\sigma\circ h)^2h_w\ov{h}_w$
(since $\Omega$ is simply connected). Therefore the
real function $ f $ is harmonic and the
map $X:=(h, f ): \Omega\rightarrow \M\times \R$ is a
conformal and minimal (possibly branched) immersion. Observe that
$X^*:=(h,- f )$ also defines a conformal and minimal (possibly branched) immersion
into $\M\times \R$, isometric to $X$ with $g^*=-g$ and
$\eta^*=-\eta$.\hfill \qed

\medskip
We denote by $\R^{2,1}$ the Minkowski 3-space, that is
$\R^3$ equipped with the Lorentzian metric $\ov{\nu}=\rmd x_1^2
+\rmd x_2^2-\rmd x_3^2$ where $(x_1,x_2,x_3)$ are the coordinates in
$\R^3$. We consider the hyperboloid $\mathcal{H}$ in $\mk$ defined by
$$
\mathcal{H}=\{(x_1,x_2,x_3)\in \mk,\ x_1^2+x_2^2 -x_3^2=-1\},
$$
$\mathcal{H}$ has two connected components, we call $\mathcal{H}_+$ the component for
which $x_3\geq 1$, we call $\mathcal{H}_-$ the other component. It is well-known that
the restriction of $\ov{\nu}$ to $\mathcal{H}_+$ is a regular metric $\nu_+$ and
that $(\mathcal{H}_+,\nu_+)$ is isometric to the hyperbolic plane $\hd$.
  We define in
the same way the metric $\nu_-$ on $\mathcal{H}_-$ and $(\mathcal{H}_-,\nu_-)$ is also
isometric to the hyperbolic plane $\hd$. Throughout this paper we
always choose as model for $\hd$ the unit disc $\di$ equipped with the
metric $\si ^2 \vert dz \vert^2=\frac{4}{(1-\vert z \vert ^2)^2}\vert dz \vert^2$.
The isometries
$\Pi_+ : \mathcal{H}_+\rightarrow \di$ and
 $\Pi_- : \mathcal{H}_-\rightarrow \di$
are given by
\begin{equation*}
\begin{split}
\Pi_+ (x_1,x_2,x_3)=&\ \frac{x_1+ix_2}{1+x_3},\ \ \forall
(x_1,x_2,x_3)\in \mathcal{H}_+ \\
\Pi_- (x_1,x_2,x_3)=&\ \frac{x_1+ix_2}{1-x_3},\ \ \forall
(x_1,x_2,x_3)\in \mathcal{H}_- \\
\end{split}
\end{equation*}
In fact $\Pi_+$ (resp. $\Pi_-$) is the stereographic projection from
the south pole $(0,0,-1)\in \mathcal{H}_-$ (resp. from the north pole
$(0,0,1)\in \mathcal{H}_+$). Observe that in \cite{[A-N]}, keeping the
notations, $\psi_2$ is the conjugate of the
stereographic projection from the south pole, that is
$\psi_2=\ov{\Pi_+}$.

Let $\Omega\subset \co$ be a connected and simply connected open
subset with
$w=u+iv$ the coordinates on $\Omega$.
An immersion $X:\Omega\rightarrow \mk$ is said to be a
{\sl spacelike immersion} if for every point $p\in \Omega$ the
restriction of $ \ov{\nu}$ at the tangent space $T_pX(\Omega)$ is a
regular metric. In this work we only consider spacelike immersions.
Let $X:\Omega\rightarrow \mk$ be a spacelike immersion.
For each $p\in \Omega$ there is a unique vector
$N(p)\in \mathcal{H}$ such that $(X_u,X_v,N)(p)$ is a positively oriented
basis and $N(p)$ is orthogonal to $X_u(p)$ and $X_v(p)$. This defines
a map $N: \Omega\rightarrow \mathcal{H}$ called the {\sl Gauss map}.

We will make now some comments about  an existence theorem of
spacelike mean curvature one surfaces in Minkowski three space  and their relation
 with harmonic maps, inferred by Akutagawa and Nishikawa [1].

Let $h:\Omega\rightarrow \di$ be a harmonic map, that is $h$ satisfies (\ref{eq:1.1}):
$$h_{w\ov{w}}+ \frac{2\ov{h}}{1-\vert h\vert^2}h_w h_{\ov{w}}=0.
$$
We assume that neither $h$ nor $\ov{h}$ is holomorphic.
It is shown in \cite{[A-N]} Theorem 6.1, that given such $h$
there exists an (possibly branched)
immersion $X_+:\Omega\rightarrow \mk$ such that the Gauss map is
$N_+=\Pi_+^{-1}\circ h$, furthermore the mean curvature is constant and
equals to 1 and the induced metric on $\Omega$ is
$$
\nu_+=\frac{4\vert h_w\vert^2}{(1-\vert h \vert^2)^2}\vert \rmd w\vert^2.
$$
We give the correspondence between our notations and the notations of
\cite{[A-N]}: $N_+ = G$ and $h=\ov{\Psi_2}$.
The map $X_+$ is unique up to a translation. In the same way there
exists an unique up to a translation (possibly branched) immersion
$X_-:\Omega\rightarrow \mk$ such that the Gauss map is
$N_-=\Pi_-^{-1}\circ h$, furthermore the mean curvature is constant and
equals to 1 and the induced metric on $\Omega$ is
$$
\nu_-=\frac{4\vert \ov{h}_w\vert^2}{(1-\vert h \vert^2)^2}\vert \rmd w\vert^2,
$$
with $N_-=G$ and $h=\Psi_1$
Let us note that those two (branched) immersions are not isometric and that
the Gauss map of $X_+$ (resp. $X_-$) takes values in
$\mathcal{H}_+$ (resp. $\mathcal{H}_-$).
In this paper we are only concerned with the immersion $X_+$.
\bigskip

On account of the above discussion we deduce:
\medskip
\begin{lem}\label{L.GaussEquation}
Let $X=(h, f ): \Omega \rightarrow \M \times \R$ be a conformal immersion. Let
$N=(N_1,N_2,N_3)$ be the Gauss map of $X$. Let $K$ (resp. $K_{ext}$)
be the intrinsic (resp. extrinsic) curvature of $X$. At last we denote
by $K_\M$ the Gauss curvature of $\M$.
Then the
{\em Gauss  equation} of $X$ reads as
$$
K(w)-K_{ext}(X(w)) = K_\M (h(w))N_3^2 (w), \leqno{(Gauss\  Equation)}
$$
for each $w\in \Omega$.
\end{lem}
\medskip
\Proof
As usual $z=x+iy$ is a local and conformal coordinate of $\M$ and $t$ is the
coordinate on $\R$. We denote by $\ov{R}$ the tensor of curvature of
$\M\times \R$, that is
$$
\ov{R}(A,B)C=\ov{ \nabla}_A \ov{\nabla}_B C -
\ov{ \nabla}_B \ov{\nabla}_A C -\ov{\nabla}_{[A,B]}C,
$$
for any vector fields $A,B,C$ on $\M\times \R$ where $\ov{\nabla}$ is
the Riemannian connection on $\M\times \R$.

As $X$ is a conformal immersion the induced metric on $\Omega$ has the
form $ds_X^2=\lambda^2 (w) \vert \rmd w^2\vert$ with
$\lambda =(\sigma\circ h)(\vert h_w\vert+ \vert  \ov{h}_w\vert)$. The
Gauss equation is
$$
K(w)-K_{ext}(X(w))=\frac{\langle \ov{R}(X_u,X_v)X_u\ ;X_v\rangle}
{\lambda^4}(w),
$$
where $\langle\ ;\ \rangle$ is the scalar product on $\M\times \R$,
$X_u=\frac{\partial X}{\partial u}=(\rme h)_u \partial_x
+(\rmi h)_u \partial_y + f _u \partial_t$ and so on. A tedious but
straightforward computation shows that
\begin{equation*}
\begin{split}
\ov{R}(\partial_x,\partial_y)\partial_x =&\
-\Delta \log (\sigma ) \partial_y\\
\ov{R}(\partial_x,\partial_y)\partial_y =&\
\Delta \log (\sigma ) \partial_x\\
\ov{R}(\partial_x,\partial_x)\partial_\ast =&\
\ov{R}(\partial_y,\partial_y)\partial_\ast=0\\
\ov{R}(\partial_t,\partial_\ast)\partial_\ast =&\
\ov{R}(\partial_\ast,\partial_t)\partial_\ast=
\ov{R}(\partial_\ast,\partial_\ast)\partial_t=0,
\end{split}
\end{equation*}
where $\partial_\ast$ stands for any vector field among
$\partial_x, \partial_y$ or $\partial_t$ and $\Delta$ is the euclidean
Laplacian. We deduce that
$$
\langle \ov{R}(X_u,X_v)X_u;X_v\rangle=
-\sigma ^2 \Delta \log(\sigma)(\vert h_w\vert^2-
\vert h_{\ov{w}}\vert^2)^2.
$$
Let us observe that $K_\M=-\Delta \log(\sigma)/\sigma^2$ therefore we
deduce from  equation (\ref{Eq.Gaussmap}) that
$N_3= (\vert h_w\vert-\vert h_{\ov{w}}\vert)/(\vert h_w\vert +
\vert h_{\ov{w}}\vert)$.
Now using the expression of $\lambda$ we get the result, which
concludes the proof.
\hfill\qed

\medskip

Notice that given a geodesic $\Gamma \subset \hi2 \times \{ 0\}$, the
vertical cylinder
$\mathcal{C}$ over $\Gamma$ defined by
$\mathcal{C} :=\{ (x,y,t);\ x+iy,\ t \in \R\} \subset \hi2 \times \R$
is a minimal surface with Gauss curvature $K \equiv 0$. We now deduce
the following.

\medskip

\begin{cor}
let $X=(h, f ): \Omega \rightarrow \hd \times \R $ be a conformal and
minimal immersion. Let $w\in \Omega$ be such that $K(w)=0$ where $K$
stands for the Gauss curvature of $X$ (that is the intrinsic
curvature).

 Then the tangent plane of $X(\Omega)$ at $X(w)$ is
 vertical. Therefore if $K\equiv 0$ then $X(\Omega)$ is part of a
vertical cylinder over a planar geodesic plane of
$\hd \times \R$, that is, there exists a
 geodesic $\Gamma$ of $\hd \times \{t \}$ such that
$X(\Omega)\subset \Gamma \times R$.
\end{cor}

\Proof
As $\M=\hd$ we have $K_\M \equiv -1$. Using the Gauss equation,
see the lemma \ref{L.GaussEquation}, we deduce that if $K(w)=0$ at
some point $w\in \Omega$ then
$$
K_{ext}(X(w))= N_3^2(w). \leqno{(\ast)}
$$
Recall that the extrinsic curvature
$K_{ext}$ is the ratio between the determinants of the
second and the first fundamental
forms of $X$. Therefore as $X$ is a minimal immersion we have
$K_{ext}(X(w))\leq 0$ at any point $w$. Using $(\ast)$
we obtain that $N_3^2(w)=0$, that is the tangent plane is vertical at
$X(w)$.

Furthermore, if $K\equiv 0$ we deduce that at each point the tangent
plane is vertical. Using this fact we get that at any point $X(w)$ the
intersection of $X(\Omega)$ with the vertical plane at $X(w)$ spanned
by $N(w)$ and $\partial_t$ is part of a vertical straight line. We
deduce that there exists a planar curve
$\Gamma \subset \hd \times \{0\}$ such that
$X(\Omega)\subset \Gamma \times \R$. Again, as $X$ is minimal we
obtain that the curvature of $\Gamma$  always vanishes, that is
$\Gamma$ is a geodesic of $\hd$. \hfill\qed

\bigskip

\section{Minimal immersions in $\M\times \R$}
Next we suppose that $\M=\m 2, \hd$ or $\sd$. In case where
$\M=\m 2$ we have $\sigma (z)\equiv 1$, if $\M =\hd$ we consider the
model of the unit disk $\di$ and then
$\sigma (z)=2/(1-\vert z\vert^2)$ for every $z \in \di$. At last if
$\M=\sd$ we can choose among the coordinate charts $\m 2$ given by the
stereographic projections with respect to the north pole and the south
pole, we have in both cases
$\sigma (z) = 2/(1+\vert z\vert^2)$ for every $z \in \m 2$.

\begin{thm}\label{T.immiso}
Let $\Omega\subset \co$ be a simply connected open set and consider
 two isometric and conformal
minimal immersions $X,X^* : \Omega\rightarrow \mathbb{M}\times \R$.
Let us call $h$ (resp. $ h^* $) the horizontal
component of $X$ (resp. $X^*$).
Assume that $h$ and $ h^* $ share the same Hopf quadratic differential.

 Then $X$ and $X^*$ are equal up to an isometry of $\mathbb{M}\times \R$.
\end{thm}

\Proof

Let us set $X=(h, f )$ where $h:\Omega\rightarrow \M$ is the horizontal
component and
$ f  : \Omega \rightarrow \R$ is the vertical component. In the same way let
us set $X^*=( h^* , f^* )$.  We will use the map $g=-ie^{\o+i\psi}$
 (resp.  $g^*=-ie^{\o^*+i\psi^*}$) associated to $h$ (resp. $h^*$) defined in the introduction
 and the one  form $\eta$ (resp. $\eta ^*$).
As $X$ and $X^*$ are isometric
immersions we infer from (\ref{eq:metric}):

$$\frac14(\vert g \vert + \vert g \vert ^{-1})^2\vert \eta \vert =
\frac14(\vert g^* \vert + \vert g^* \vert ^{-1})^2\vert \eta^* \vert$$
Also as $h$ and $ h^* $ share the same Hopf quadratic differential $Q=\phi dw^2$ we have
$$\vert \eta \vert =2 \vert \phi \vert ^{1/2}\vert= \vert \eta ^* \vert $$
We deduce that we have
$$
\vert g \vert = \vert g^* \vert ,\leqno{(\ast})
$$
or
$$
\vert g \vert = \vert g^* \vert ^{-1},\leqno{(\ast\ast)}
$$
If case $(\ast\ast)$ happens we consider the new immersion
$X^{**} : \Omega\rightarrow \M \times \R$ defined by $X^{**}=(\ov{ h^* }, f^* )$. Now case
$(\ast)$ happens considering immersion $X^{**}$, with datas
$g^{**}= (g^{*})^{-1}$ and $\eta ^{**}=\eta ^{*}$.
Note that $X^{**}$ and
$X$ are isometric immersions with same Hopf quadratic differential.
Therefore, up to an isometry of
$\M \times \R$, we can assume that case $(\ast)$ happens
and $\omega = \omega ^*$.

\medskip

Let us assume now that $\M=\hd$, the case $\M=\sd$ is similar and case
$\M=\m2$ will be considered later.

Let us consider the Minkowski 3-space $\R^{2,1}$.
 As $h:\Omega \rightarrow \hd$ is a harmonic map and $\Omega$ is simply connected
 it is known that there exists a CMC one
(possibly branched) immersion
$\wt{X}: \Omega\rightarrow \R^{2,1}$ such that the Gauss map is
$\Pi_+^{-1} \circ h$.
Furthermore the induced metric on $ \Omega $ is given by
$$
ds^2_{\wt{X}}=\left( (\sigma \circ h) \vert h_w\vert\right)^2
\vert \rmd w\vert^2= e^{2 \omega}\vert \phi \vert\vert \rmd w\vert^2 ,
$$
see the Section \ref{Sec.Prel} (Preliminar).
Notice that $\phi$ can vanish only at isolated points, so there
exists a simply connected open subset $V$ of $\Omega$, $V\subset \Omega$, such
that $\wt{X}$ defines a regular immersion from $V$ into $ \R^{2,1}$ and
$ds^2_{\wt{X}}$ defines a regular metric.

 Furthermore we deduce from Theorem 3.4 of \cite{[A-N]} that the
 second fundamental form of $\wt{X}$ is given uniquely in term of $Q$ and
$ds^2_{\wt{X}}$. To see this, observe first that, setting
$\tilde{ \phi}(\wt{X}):=\frac12(b_{uu}-b_{vv} -i2 b_{uv})$, we get from relation (3.12) of
\cite{[A-N]} that
\begin{equation}\label{Eq.Jap}
\tilde{\phi}(\wt{X})=(\sigma\circ h)^2 h_w \ov{h}_w=\phi,
\end{equation}
that is $\tilde{ \phi}(\wt{X})\rmd w^2=Q(h)$. Pay attention to
the fact that our
conventions are not the same as in \cite{[A-N]}, for example following
notation of \cite{[A-N]} we have \newline
$b_{uu}= \left( (\sigma \circ h) \vert h_w\vert\right)^2 h_{11}$ and
so on, therefore
$\tilde{ \phi}(\wt{X}) =\left( (\sigma \circ h) \vert h_w\vert\right)^2
\phi$ where $\phi $ is given in \cite{[A-N]} that is
$\phi =\frac12(h_{11}-h_{22} -i2h_{12})$.

 Now using the fact that
$b_{uu}+b_{vv}=2\left( (\sigma \circ h) \vert {h}_w\vert\right)^2$
(since the mean curvature is 1), we deduce that
\begin{eqnarray}
b_{uu}(\wt{X}) &=& \left( (\sigma \circ h) \vert {h}_w\vert\right)^2
+ \rme \ Q(h)/\rmd w^2=e^{2\o}\vert \phi \vert +\rme \phi  \label{Eq1} \\
b_{vv}(\wt{X}) &=& \left( (\sigma \circ h) \vert {h}_w\vert\right)^2
- \rme \ Q(h)/\rmd w^2= e^{2\o}\vert \phi \vert -\rme \phi  \label{Eq2}  \\
b_{uv}(\wt{X}) &=& -\rmi \ Q(h)/\rmd w^2=-\rmi \phi \label{Eq3}
\end{eqnarray}

In the same way there exists an unique (up to a translation) CMC one
(possibly branched) immersion $\wt{X^*} :\Omega\rightarrow \R^{2,1}$ such
that the Gauss map is
$\Pi_+^{-1} \circ  h^* $. We can assume that $\wt{X^*}$ defines a regular immersion on
$V$. Notice that we have $Q (h)=Q ( h^* )$ and identities
$(\ast)$ as well. We deduce
from the former discussion that
$\wt{X}$ and $\wt{X^*}$ share the same induced metric on $V$ and the
same second fundamental form. Therefore we infer with the fundamental
theorem of geometry in Minkowski 3-space that $\wt{X}$ and
$\wt{X^*}$ are equal up to a
positive isometry $\Gamma$ in $\R^{2,1}$, that is
$\wt{X^*}=\Gamma \circ \wt{X}$. The restriction of $\Gamma$
on $\hd$ defines an isometry $\gamma$ of $\hd$ and we get
$ h^* =\gamma \circ h$ on $V$. By an argument of analyticity we have
$ h^* =\gamma \circ h$ on the entire $\Omega$.

 Let us return to $\hd \times \R$. As $ f^* _w=\pm  f _w$ in view of
(\ref{eq:third})
 we get that
$ f^* =\pm  f  +c$ where $c$ is a real constant. At last we obtain
$X^*:=( h^* , f^* )=(\gamma \circ h, \pm  f  +c)$, that is $X^*$ and $X$ differ from
an isometry of $\hd \times \R$.

\medskip

In the case where $\M=\sd$ the proof is similar: we use the fact
that any harmonic map from $\Omega$ into $\sd$ is the Gauss map of an
unique (up to a translation) CMC 1 (possibly branched) immersion into
$\m3$, see \cite{[Ken]}.

\medskip

Finally let us consider the case where $\M=\m2$. Let $(g,\eta )$
(resp. $(g^*,\eta^* )$) be the Weierstrass representation of
$X$ (resp. $X^*$). Therefore $X$ is given by
$X= \left(
\frac12\ov{\int  g^{-1}\eta}-\frac12\int g \eta ,
 \rme \int \eta \right)$.
As $\vert g^*\vert =\vert g\vert$, we deduce that there exists a real
number $\theta$ such that $g^*=e^{i\theta}g$. Furthermore we have
$\eta= \pm \eta^*$ since $(f_{z})^2=(f^*_{z})^2=-\phi$. Finally we have
$(g^*,\eta^*)=(e^{i\theta}g,\pm \eta )$ and we
deduce that $X^*$ differ from $X$ by an isometry of $\m2\times\R$,
this concludes the proof. \hfill\qed

\bigskip

There is also an existence result of minimal immersion into
$\M \times \R$ where $\M= \hd, \sd$ or $\m 2$.

\begin{thm}\label{T.existence}
Let $\Omega \subset \co$ be a simply connected domain. Let
$ds^2=\lambda^2(w)\vert \rmd w\vert^2$ be a conformal metric on $\Omega$
and let $Q =\phi (w)\rmd w²$ be a holomorphic quadratic differential on
$\Omega$ with zeros (if any) of even order. Assume that $\M= \hd, \sd$
or $\m 2$.

  Then there exists a conformal and minimal immersion
$X:  \Omega \rightarrow  \M \times \R$  such that, setting
$X:=(h, f )$, the Hopf
  quadratic form of $h$ is $Q$ (that is $Q (h) =\phi (w)\rmd w^2$)
  and such that the induced metric $ds^2_X$ is
$$
ds^2_X =ds^2=\lambda^2(w)\vert \rmd w\vert^2
$$
if and only if $\lambda$ satisfies
$\lambda^2 -4 \vert \phi \vert \geq 0$ and
\begin{equation}
\label{Eq.cond}
\Delta \o =-2K_{\M} \sh 2\o \vert \phi \vert
\end{equation}
where $K_{\M}$ is the (constant) Gauss curvature of $\M$ and
$$ \o :=\log \frac{\lambda  -\sqrt{\lambda^2-4\vert \phi\vert}}{2}
-\frac12 \log \vert \phi \vert$$
\end{thm}

\Proof
We first consider the case $K_{\M}=-1$, that is $\M=\hd$.
Let us assume that $\lambda$ satisfies (\ref{Eq.cond}). Consider the
2-form $II:= b_{uu}du^2+2b_{uv}dudv+b_{vv}dv^2$ on $\Omega$ where
$b_{uu},b_{uv}$ and $b_{vv}$ are given by:
\begin{equation}\label{Eq.second}
\left\{ \begin{split}
b_{uu} +b_{vv}=&\ 2 e^{2\o}\vert \phi \vert \\
b_{uu} -b_{vv}=&\ 2 \rme (\phi)\\
b_{uv}=&\ -\rmi (\phi)
\end{split}\right.
\end{equation}
The Gauss equation for the pair $(e^{2\o}\vert \phi \vert \vert \rmd w\vert^2, II)$
in $\R^{2,1}$ is:
$$\Delta \o =-2 \sh (2\o) \vert \phi \vert$$
and then it is satisfied.
The Codazzi-Mainardi equations are also satisfied since $\phi$ is
holomorphic. Therefore the fundamental theorem of geometry in
$\R^{2,1}$ states that there exists an immersion
$\wt X :\Omega \rightarrow \R^{2,1}$ such that the induced metric on
$\Omega$ is $ds^2_{\wt X}=e^{2\o}\vert \phi \vert \vert \rmd w\vert^2$
and the second fundamental form is $II$. Now the equations (\ref{Eq.second}) show
that the immersion has constant mean curvature one.

 Up to an isometry of $\R^{2,1}$ we can assume that the Gauss map $N$
 of $\wt X$ takes values in $\mathcal{H}_+$. Therefore
$h:=\Pi_+ \circ N:  \Omega\rightarrow \hd$ is a harmonic mapping  such
  that its Hopf quadratic form is the same as $\wt X$:
$Q (h)=\tilde{\phi}  (\wt X)\rmd w^2$, as we have seen in the proof of
Theorem \ref{T.immiso}, see relation (\ref{Eq.Jap}). By definition we have:
$$
\tilde{Q} (\wt X)\rmd w^2:=\frac12(b_{uu}-b_{vv}-i2b_{uv})\rmd w^2=Q.
$$
Therefore we obtain $Q (h) =Q$. Moreover we have
$$
ds^2_{\wt X}=\left( (\sigma \circ h) \vert h_w\vert\right)^2
\vert \rmd w\vert^2
$$
and we deduce that
$e^{2\o}\vert \phi \vert =\left( (\sigma \circ h) \vert h_w\vert\right)^2$.

Now we apply Proposition
\ref{P.min.imm}  which states that there exists a conformal and minimal
immersion $X=(h, f ):\Omega \rightarrow \hd \times \R$, with induced
metric:
$$
ds^2_X=(\sigma^2 \circ h) (\vert h_w\vert + \vert \ov{h}_w\vert)^2
\vert \rmd w\vert^2
$$
At last using the fact that
$(\sigma \circ h) \vert \ov{h}_w\vert=
\vert \phi\vert/(\sigma \circ h) \vert h_w\vert$ we easily compute
that:
$$
(\sigma^2 \circ h) (\vert h_w\vert + \vert \ov{h}_w\vert)^2=
\ch ^2 \o \vert \phi \vert \vert \rmd w\vert^2=\lambda^2,
$$
that is $ds^2_X=\lambda^2 \vert \rmd w\vert^2$ as desired.

 Conversely suppose that such an immersion exists. Then we have by (see (\ref{eq:metric})):
$$
\lambda ^2 =4\ch^2 \o \vert \phi \vert \vert \rmd w\vert^2 .
$$
A simple computation shows that we have:
$$
\o =\o_1:=\frac12\log\frac{\vert \ov{h}_w\vert}{\vert h_{w}\vert} \ \
\mathrm{or}\ \
\o = \o_2:=\frac12\log \frac {\vert h_w\vert}{\vert \ov{h}_{w}\vert}.
$$

The equation  (\ref{Eq.cond}) is B\"ochner formula (\ref{eq:bochner}).
This concludes the proof in case where $\M=\hd$.

  If $\M=\sd$ (and then $K_{\M}=1$),  the proof is analogous: we
  use the fact that for any constant mean curvature one immersion
$\wt X :\Omega \rightarrow \m 3$ its Gauss map
$N:\Omega \rightarrow \sd$ is harmonic and conversely any harmonic map
from $\Omega$ into $\sd$ is the Gauss map of an (possibly branched)
  immersion into $\m 3$ with constant mean curvature one.

  If $K_{\M}=0$, that is $\M =\m 2$, assume first that $\o$ satisfies
  (\ref{Eq.cond}), that is $\o$ is a harmonic function. As
  $\Omega$ is simply connected, $\o$ is the real part of a
  holomorphic function $\o+i\psi$ on $\Omega$. We set:
$$\eta:=-2i\sqrt \phi\ \ \mathrm{and}\ \ g:=-ie^{\o+i\psi}.$$
Let $X=(h, f ):\Omega \rightarrow \m 2 \times \R$ be the conformal and
minimal immersion given by the Weierstrass representation
$(g,\eta)$.
\hfill\qed

\medskip

Observe that in case where $K_{\M}=0$ the result can be encountered in
\cite{[Du]}, see the theorem in Section 10.2. We gave the proof for
sake of completness. The cases $K_{\M}=1$ and $K_{\M}=-1$ were proved by
B. Daniel in \cite{[Dan]} using other methods.

\bigskip

\begin{defn}\label{D.asso}
Let $M$ be any riemannian surface.
Let $X,X^* : \Omega\rightarrow \M \times \R$ be two conformal minimal
immersions and let us set
$X=(h, f )$ and $X=( h^* , f^* )$.

For any $\theta \in \R$ we say that $X$ and $X^*$ are
{\sl$\theta$- associate} (or simply {\sl associate}) if
they are isometric
immersions and if  we have
$Q( h^* )=e^{2i\theta}Q(h)$. That is $X$ and $X^*$ are associate if
and only if we have
$$
(\sigma \circ h) (\vert h_w\vert + \vert \ov{h}_w\vert)=
(\sigma \circ  h^* ) (\vert  h^* _w\vert + \vert \ov{ h^* }_w\vert)\ \
\mathrm{and}\ \
(\sigma\circ  h^* )^2  h^* _w \ov{ h^* }_w=
e^{2i\theta}(\sigma\circ h)^2 h_w \ov{h}_w,
$$
where, in a local coordinate $(z)$, the metric on $\M$ is given by
$\sigma ^2(z) \vert \rmd z \vert ^2$. In case where
$\M=\m 2, \hd$ or $\sd$, we
deduce from the theorem \ref{T.immiso} that given a conformal
minimal immersion $X$, the $\theta$-associate minimal immersion is
uniquely determined up to an isometry of $\M\times\R$.
Furthermore if $\theta=\pi/2$ we say that $X$ and $X^*$ are
{\sl  conjugate}.
\end{defn}

\begin{rem}
\label{associate} Two isometric immersions $X$ and $X^{\te}$ are associate up to an isometry if
$\eta^{\te}=e^{i\theta}\eta$ and by (\ref{eq:metric})
$\vert g^{\te} \vert + \vert g^{\te} \vert ^{-1}=
\vert g \vert + \vert g \vert ^{-1}$
(or equivalently $\ch \o^{\te}=\ch \o$). Then $\o^{\te}=\o$ or $\o^{\te}=-\o$.
In particular $X$ and $X^{\te}$ are associate if and only if $N_{3}(X)=N_{3}(X^*)$
or $N_{3}(X)=-N_{3}(X^*)$ (recall that $N_{3}(X)=\ta \o$) and
$\eta ^{\te}= e^{i \te} \eta$.
\end{rem}

\bigskip

 In fact B. Daniel proved that the associate family always
exists in $\hd \times \R$ and $\sd \times \R$, see \cite{[Dan]}. In
this situation he gave an alternative definition of associate and
conjugate  isometric immersions which turn to be equivalent to our definition.
We are going to give another proof of the existence of the associate
family.

\bigskip

\begin{cor}\label{C.associate}
Let $X:=(h, f ):\Omega\rightarrow \M \times \R$ be any conformal and
minimal immersion where $M=\hd,\ \sd$ or $\m 2$. Then for any
$\theta\in \R$ there exists a $\theta-associate$ immersion
$X_\theta :=(h_\theta, f _\theta):\Omega\rightarrow \M \times \R$.
Furthermore $X_0=X$ and $X_\theta$ is unique up to isometry of
$\M\times \R$.
\end{cor}
\Proof
Let us set $Q (h)=\phi (w)dw^2$ and let $ds^2_X$ be the
conformal metric induced on $\Omega$  by $X$. We deduce from
Theorem \ref{T.existence} that the pair $(ds^2_X, \phi)$ satisfies the
condition (\ref{Eq.cond}). Therefore for any $\theta \in \R$ the pair
$(ds^2_X,e^{2i\theta} \phi)$ also satisfies condition
(\ref{Eq.cond}). Finally we infer with Theorem \ref{T.existence} that
there exists a $\theta$-associate immersion, which concludes the
proof. \hfill\qed

\section{Minimal vertical graph}
In this section we study geometric properties of minimal graph and their associate
family. Recall from the introduction that we introduce some "Weierstrass" data
for minimal surfaces $(g,\eta)$ with $g=-ie^{\o+i\psi}$ and $\eta=-2i \sqrt \phi$.
When $X$ is a minimal surface of $\R ^3$, then $\o+i \psi$ is meromorphic. In
the other case $\o$ satisfy the $\sh$-Gordon equation (\ref{eq:gordon}). In  the following
Lemma, we determine how the function $\o+i\psi$ deviate from to be meromorphic.
We express same expression for associate family. In this case (remark \ref{associate}),
up to an isometry we have $g^{\te}=-ie^{\o+i\psi^{\te}}$ and
$\eta ^{\te}=e^{i\te}\eta$ ($\o^{\te}=\o$).
Then we have:

\begin{lem}
\label{cauchy}
We consider  a harmonic map $h :\Omega \rightarrow (U,\si ^2 \vert \rmd z \vert ^2)$
with  holomorphic quadratic Hopf differential $Q=\phi(w) (\rmd w)^2$ with zeros
(if any) of even order and coefficient of dilatation $a(z)=e^{-2(\o+i\psi)}$.
Thus  we can define $\sqrt \phi =\vert \phi \vert^{1/2} e^{i\beta}$ and we identify
 $\si$ with $\si \circ h$. Then
\begin{equation}
(\o+i\psi)_{\bar w}=
\vert \phi\vert ^{1/2}e^{-i\be}\left( \sh \o
\langle \frac{\nabla  \log \si}{ \si}, e^{i\psi}\rangle
 +i \ch \o\langle \frac{\nabla  \log \si}{ \si},ie^{i\psi}\rangle \right)
\label{eq:cauchy1}
\end{equation}
\end{lem}
\begin{cor} If $X=(h,f)$ is a minimal surface and
$X^{\theta}=(h^{\theta},\eta^{\theta})$ is the associate family of $X$ define in the
Definition \ref{D.asso}, we can define the map $\o^{\te}+i\psi^{\te}$ and with the
notation $\si^{\te}=\si \circ h^{\te}$ we have $\o^{\te}=\o$
and

\begin{equation*}
(\o+i\psi^{\te})_{\bar w} =
\vert \phi\vert ^{1/2}e^{-i(\be+\te)}\left( \sh \o
\langle \frac{\nabla  \log \si^{\te}}{ \si^{\te}}, e^{i\psi^{\te}}
\rangle + i \ch \o \langle \frac{\nabla  \log \si^{\te}}{ \si^{\te}},
ie^{i\psi^{\te}}\rangle
\right)
\end{equation*}

\end{cor}

\Proof We compute $\psi _u$ as a  function of $\o_v$ and $\psi _v$
as a function of $\o_u$. In complex coordinates $w$, using
(\ref{eq:1}), (\ref{eq:gauss}) and assuming $\eta =-2i
\sqrt{\phi}$ we derive:

$$\f _w = \frac{\sqrt \phi e^{\o +i\psi}}{ \si}
 \hbox{ and } \f _{\bar w}=
 \frac{\overline{\sqrt \phi} e^{-\o +i\psi}}{ \si} $$
 \noindent
 while
 $$ \f ^{\te}_{w} =
\frac{ e^{i\te} \sqrt \phi e^{\o +i \psi^{\te}}}{ \si^{\te}}
 \hbox{ and } \f^{\te} _{\bar w}=
 \frac{e^{-i\te}\overline{\sqrt \phi} e^{-\o +i\psi^{\te}}}{ \si^{\te}}.$$
\noindent
Inserting these expressions in the harmonic equation $(\ref{eq:1.1})$ we obtain:

$$(\o + i \psi )_{\bar w} = -\si\left(\frac{1}{\si}
\right)_{\bar z} - 2(\log \si )_u \f _{\bar w}$$

$$(\o + i \psi ^{\te})_{\bar w} = -\si^{\te}\left(\frac{1}{\si^{\te}}
\right)_{\bar w} -2 (\log \si^{\te} )_u \f _{\bar w}.$$
\noindent
Now note that

\[
\begin{array}{lll}
\displaystyle{ - \si\left(\frac{1}{ \si} \right)_{\bar w} }&
=&\displaystyle{ (\log\si)_{\bar w} }\cr
 & = & \displaystyle{\left( (\log \si)_z \f _{\bar w} + (\log
\si)_{\bar z}{\bar \f } _ {\bar w} \right) }
\end{array}
\]
\noindent
where $2(\log \si)_{z}=(\log \si)_{x} -i(\log
\si)_{y}$ and $ \displaystyle { {\bar \f } _ { \bar z}= \overline{\f _{z}}}$.
 Collecting these equations we obtain:

$$(\o + i \psi )_{\bar w} =
 (\log \si)_{\bar z} {\bar\f}_{\bar w}-
 (\log \si )_z \f_ {\bar w}. $$

\noindent
which is

\[
\begin{aligned}
(\o + i \psi )_{\bar w} &=\frac { \vert \phi\vert^{1/2} e^{-i\beta}}
 {  \si}\left(
 \sh \o \left(\cos \psi  (\log \si)_{x}  + \sin \psi(\log \si )_{y} \right)\right.
 \\
 &\quad \left. +i  \ch \o
 \left(\cos \psi  (\log \si)_{y}  -\sin \psi(\log \si )_{x} \right)\right)
\end{aligned}
\]
\noindent
Since $X^{\te}$ is isometric to $X$, we have $\o^{\te}=\o$
and the same equation applied to $\f^{\te}$ yields,

$$(\o + i \psi ^{\te})_{\bar w} =
  (\log \si^{\te})_{\bar z} {\bar\f^{\te}}_{\bar w}-
 (\log \si ^{\te})_z \f^{\te}_ {\bar w}. $$

\noindent
gives
\[
\begin{aligned}
(\o + i \psi ^{\te} )_{\bar w} &=\frac { \vert \phi\vert^{1/2} e^{-i(\beta+\te)}}
 {  \si^{\te}}\left(
 \sh \o
 \left(\cos \psi^{\te}  (\log \si^{\te})_{y}  +
\sin \psi^{\te}(\log \si ^{\te})_{y}
  \right)\right.
 \\
 &\quad \left. -i
 \ch \o \left(\cos \psi^{\te}  (\log \si^{\te})_{y}  -
  \sin \psi^{\te}(\log \si^{\te} )_{y}
   \right)\right)
\end{aligned}
\] \hfill\qed

We consider the projection
$F :\M \times \R \longrightarrow \M \times \{ 0\}$, thus
$F \circ X= \f$. Now let us consider a curve $\gamma :[0,l]
\longrightarrow \Om \subset \C$ parametrized by arclenght and
$\ga '(t)=e^{i\alpha (t)}$ in $\Om \subset  \C$. We will
compute in the following what are the curvature $k$  in $\M$ of
the planar curves $F \circ X (\gamma)=h(\gamma)$ and $F \circ
X^{\te}(\gamma)=h^{\te}(\ga)$. A such computation appears in \cite{[H]} in
the particular case where $\alpha = 0$ and $\alpha = \pi/2$:
\begin{prop}
\label{courbure}
If we consider a curve $\ga$ in $\Om$ and the image
 $\f(\ga))$ and $\f^{\te}(\ga)$ in $M$, then
the curvature are

\begin{equation}
\label{eq:courbure}
k(\f(\ga))  =
\frac{\sin \al \o_{u}-\cos \al \o _{v} + G_{t}}{2\vert \phi \vert ^{1/2}R}
\end{equation}

\begin{equation}
\label{eq:courbureteta}
k(\f^{\te}(\ga)) =
\frac{\sin \al \o_{u}-\cos \al \o _{v} + G^{\te}_{t}}{2\vert \phi \vert ^{1/2}R^{\te}}
\end{equation}

 \noindent
 where
\[
\begin{array}{l}
Re^{iG}=\cos (\al+\be)\ch \o
+ i \sin (\al+\be) \sh \o \\
\\
R^{\te}e^{iG^{\te}}=\cos (\al+\be+\te)\ch \o
+ i \sin (\al+\be+\te) \sh \o
\end{array}
\]
\end{prop}
\Proof
We apply the formula (\ref{eq:1}) with $g=-ie^{\o+i\psi}$, $\eta
=-2i\sqrt {\phi} dz$ for $X$ and $g^{\te}=-ie^{\o+i\psi^{\te}}$
and $\eta ^{\te}=e^{i\te}\eta=-2ie^{i\te}\sqrt{\phi}dz$. Let us recall
 $\sqrt{\phi}=\vert \phi\vert ^{1/2}e^{i\beta}$.

\[
\begin{aligned}
\frac{d h(\ga)}{dt} & =
\frac{2\vert \phi \vert ^{1/2}}{ \si}
\ch (\o + i \alpha + i \beta) e^{i\psi} \\
& = \frac{2\vert \phi \vert ^{1/2}}{ \si}(\cos (\al+\be)\ch \o
+ i \sin (\al+\be) \sh \o) e^{i\psi} \\
\frac{ d h(\ga)}{dt} & = \frac{2 \vert \phi \vert ^{1/2}}{\si}
R e^{i(\psi + G)}
\end{aligned}
\]

\[
\begin{aligned}
\frac{d h^{\te}(\ga)}{dt} & =  \frac{2 \vert \phi \vert^{1/2}}{\si}
\ch (\o +i\beta+i\al+i\te)e^{i\psi^{\te}}\\
& = \frac{2  \vert \phi \vert^{1/2}}{ \si} (\cos (\al+\be+\te)\ch \o
+ i \sin (\al+\be+\te) \sh \o) e^{i\psi^{\te}}\\
\frac{ d h^{\te} }{d t} & = \frac{ 2 \vert \phi \vert ^{1/2}}
{\si} R e^{i(\psi^{\te} + G_{\te})}
\end{aligned}
\]


\noindent


If $k$ is the curvature of a curve $\gamma $  in $(U, \si^2 (z) \vert \rmd z
\vert^2)$ and $k_{e}$ is the Euclidean curvature in
$(U,\vert dz \vert^2)$,
we get by conformal change of the metric:

$$k= \frac{k_e}{\si}-\frac{ \left< \nabla \si , n \right>
}{\si^2}$$
\noindent
where $n$ is the Euclidean normal to the curve $\gamma $ such that $(\ga ', n)$ is
positively oriented. If
$s$ denotes the arclength of $\f(\gamma) $ and $s^{\te}$ the
arclength of $\f^{\te}(\ga)$, we have

\begin{equation}
\label{eq:courb.eucli}
k_e(\f(\ga) )= \psi _s+G_{s} = \frac{\si}{2 \vert \phi \vert ^{1/2} R} (\cos
\al \psi _u+\sin \al \psi_{v})+G_{s}
\end{equation}

\noindent
We consider the euclidean normal of $\f(\ga)$ (resp  $\f^{\te}(\ga)$) given
by

 \[
 \begin{array}{l}
\displaystyle{n=(-\sin (\al +\be)\sh \o + i \cos (\al +\be)\ch \o)\frac{e^{i\psi}}{R}} \\
\displaystyle{n^{\te}=(-\sin (\al +\be+\te)\sh \o+ i \cos (\al +\be+\te)\ch \o)
\frac{e^{i\psi^{\te}}}{R^{\te}}}
\end{array}
\]
\noindent
and
 \[
 \begin{aligned}
 \frac{ \left< \nabla \si , n \right>}{\si^2}=\frac{ \left< \nabla
\log \si , n
\right>}{ \si}=-\sin (\al+\be)\frac{\sh \o}{R} \langle
\frac{\nabla \log \si}{ \si}, e^{i\psi}\rangle \cr
+\cos (\al+\be)\frac{\ch \o}{R}
 \langle \frac{\nabla \log \si}{ \si}, i e^{i\psi}\rangle
 \end{aligned}
 \]
Using the Lemma \ref{cauchy}, we express $\psi _u$ in terms of
$\o_v,$ and $\psi _v$ in terms of $\o _u$ by using formula
(\ref{eq:cauchy1}) in (\ref{eq:courb.eucli}), we deduce

$$\frac{\psi _{s}}{ \si} -\frac{ \left< \nabla \si,
 n \right>}{\si^2}  = \frac{\sin \al \o_{u}-\cos \al \o _{v}}
 {2\vert \phi \vert ^{1/2}R} $$

The same computation with $X^{\te}$ yields
$$\frac{\psi ^{\te}_{s^{\te}}}{\si^{\te}} -\frac{ \left< \nabla {\si^{\te}},
 n \right>}{(\si^{\te})^2}  = \frac{\sin \al \o_{u}-\cos \al \o _{v}}
 {2\vert \phi \vert ^{1/2}R^{\te}} $$

This concludes the proof of the proposition since $G_{s}=\frac{\si}{2 \vert \phi \vert ^{1/2} R}G_{t}$ . \hfill\qed

 Now we prove the generalization of the Krust's theorem for minimal vertical graph and  associate  family
 surfaces. Let $U\subset \M$ be an open set and $f(z)$ a smooth function on $U.$
We say that $G$ is a vertical graph in $\M \times \R$ if
$G=\left\{ (z,t) \in \M \times \R ; t=f(z),\, z\in U \right\}$. The graph
is an entire vertical graph if $U=\M$. 

\begin{thm}
\label{krust}
If we consider a minimal graph $X (\Om)$ on a convex domain $\f(\Om)$ in $\M$,
then the associate surface $X^{\te}(\Om)$ is a graph under the assumption that the
curvature $K_{\M}\leq 0$.
\end{thm}
When $K_{\M}\equiv 0$ it is a result of Romain Krust
(see \cite{hil}, page 188 and application therein).

\Proof
The proof is a direct application of Gauss-Bonnet theorem with the fact that
$\o$ is without zero ($X$ is a vertical graph). If we consider a smooth
piece of  embedded curve $\Gamma$ in $\M$ with end points $p_{1}$ and $p_{2}$, then if
$p_{1}=p_{2}$, $\Gamma$ is enclosing a region $A$ and:

$$\int_{A}K_{\M}dV_{\si}+\int_{\Gamma}k(s) ds + \al=2\pi$$

\noindent
where $\al$ is the exterior angle  at $p_{1}=p_{2}$. Since
$\Gamma$ is embedded we have
$\al \leq \pi$. The Gauss Bonnet formula above gives us that in the  case
where $K_{\M}\leq 0$:

$$\pi\geq \al \geq 2\pi-\int_{\Gamma}k ds.$$

Now, if we assume that $X^{\te}(\Om)$ is not a graph, there exist two points $p_{1}$ and $p_{2}$
two points of $\Om$ with $\f ^{\te}(p_{1})=\f^{\te}(p_{2})$. Since $\f (\Om)$
is convex, there is a geodesic in $\f(\Om)$ which can be lift by a path $\ga$ in $\Om$.
In summary, we assume that the curve $\ga(t)$, $t \in [0,l]$
is parametrized by euclidean arclength, $\ga '(t)=e^{i\al}$,
$\f(\ga)$ is a piece of a geodesic of $\M$, $p_{1}$ and $p_{2}$ are the end points of $\ga$
and $\f^{\te}(p_{1})=\f^{\te}(p_{2})$. We assume that $\f^{\te}(\ga)$
is a closed embedded curve. If not, we can consider a subarc of $\gamma$ with end
points $p'_{1}$ and $p'_{2}$,
with an image by $\f^{\te}$  smooth, embedded and $\f^{\te}(p'_{1})=
\f^{\te}(p'_{2})$. In
the case where this subarc embedded doesn't exist, then it is meaning that all points are
double, like a path where we go and back after an interior end point $q$. At $q$,
$\f^{\te}(\ga)$ is not immersed. The derivative is nul and then the tangent plane of
$X^{\te}$ is
vertical. Then $\o$ would have an interior zero (a contradiction with the vertical graph assumption).

We will apply the formula of Gauss-Bonnet and we will prove that
$\int_{\f^{\te}(\ga)} k ds^* < \pi$ under the  hypothesis that
$\f(\ga)$ is  a geodesic. It will provide  a contradiction with
$\al \leq \pi$ and then the horizontal curve cannot be  closed and
embedded.

If $\f(\ga)$ is a geodesic, then by the formula
(\ref{eq:courbure}) of the previous proposition \ref{courbure}, we have
$\sin \al \o_{u}-\cos \al \o _{v} + G_{t}=0$. Thus

$$k(\f^{\te}(\ga))=  \frac{ G^{\te}_{t}-G_{t}}
 {2\vert \phi \vert ^{1/2}R^{\te}}$$

 Since $ds^*=2\vert \phi \vert ^{1/2}R^{\te} dt $, we have
 $$\int_{\f^{\te}(\ga)} k ds^*=( G^{\te}(l)-G(l))-(G^{\te}(0)-G(0))$$

 Now, we remark with a direct computation of the real and imaginary part
 of

 $$z=\frac{R^{\te}}{R}e^{i(G^{\te}-G)}=\frac{\cos(\al+\be+\te) \ch \o +
 i \sin(\al+\be+\te) \sh \o}
 {\cos (\al+\be)\ch \o + i \sin (\al+\be)\sh \o}$$
 \noindent
 that
 $$\tan(G^{\te}(t)-G(t))=\frac{\sh (2\o)\sin \te}{2 \cos \te R^2 -
 \sin \te \sin 2 (\al+\be)}$$.

 Since $X$ is a graph $\o$ is  without interior zero, and then
 ${\sh (2\o)\sin \te}$ cannot be zero for $\te \in (0,\pi/2]$.
 It implies that $G^{\te}(t)-G(t) \in [0, \pi]$. \hfill\qed

The example 13 below prove that the conjugate surface of an entire graph in $\hd \times \R$
is a graph (by our theorem \ref{krust}) but it is not necessary entire. In this
direction we give the following criterion in $\hd \times \R$:

\begin{thm}
\label{entier} 
Let $X: \di ^{2} \longrightarrow \hd \times \R$ be an entire vertical graph on
$\hd$. Then  on any divergent path $\ga$ of finite euclidean lenght in  $\di$

$$\int_{\ga} \vert \phi  \vert ^{1/2} dt < \infty$$
\noindent
then the conjugate graph $X^*$ is an entire graph.
\end{thm}

\Proof
Recall that $f+if^*=-2i\int ^{z}\sqrt{\phi}$ is holomorphic. 
We consider $\ga (t)$ a divergent
path in $\di ^{2}$ and $X(\ga)=\Gamma$ its image in the graph (recall that $\ga ' (t)=e^{i\al})$). Since $X$ is a proper
map, the length of $\Gamma$ is infinite in $X$ and

\begin{equation}
\label{long}
\ell (\Gamma)=\int_{\ga} 2 \cosh \o \vert \phi \vert ^{1/2} dt =\infty
\end{equation}

Now we prove that  the length of $h^*\circ \ga$ is infinite which prove the theorem.
If $X^*$ is not entire, one can find a diverging curve in 
$\di ^2$ with $h^* \circ \ga$ of finite length. Now let us compute

\[
\ell (h ^* \circ \ga )=\int _{\ga } 2 \vert  \phi \vert ^{1/2} R^{*} dt
\]

\noindent
where $R^{*2}=\sin ^2 (\al + \be) \cosh ^2 \o + \cos ^2 (\al + \be) \sinh ^2 \o$
(recall  $R^{*}=R^{\pi/2}$). Now we remark that

$$R^{*2} = \cosh ^2 \o - \cos ^2 (\al + \be).$$
\noindent
Then, using (\ref{long}), and the hypothesis we have

$$\ell (h ^* \circ \ga )= \ell(\Gamma) -
\int_{\ga}\vert  \cos (\al+ \be)\vert  \vert \phi \vert ^{1/2} dt =\infty .$$ 

\hfill\qed
\section{Examples}

\begin{exa}\label{Ex.Scherk}
Let us consider the {\em Scherk type surface} in
$\hd \times \R$
 invariant by hyperbolic translations. It is a complete minimal graph over a
non-bounded domain in $\hd$ defined by a complete geodesic $\gamma$ in
$\hd \times \{0\}$. The graph takes $\pm \infty$ value on $\gamma$ and $0$
 value on the asymptotic boundary. In the half plane model of
$\hd=\{(x,y)\in \m 2, \ y>0\}$, there is a nice formula for the graph:

\begin{equation*}\label{sch}
t = \ln \left( \frac{\sqrt{x^2 +y^2} + y}{x}\right),\qquad y>0,\
x>0
\end{equation*}

 The  Scherk's conjugate minimal
 surface in $\hd \times \R$ (see Theorem 4.2 of \cite{[Earp]})
 is given by the equation: $t=x$.
It is invariant by parabolic screw-motions. It is a whole graph over $\hd$.
The second and third authors proved
that in $\hd\times \R,$ a catenoid is conjugate to a
helicoid of pitch $\ell <1$, see \cite{[S-T]}. Surprisingly,
a helicoid of pitch
$\ell =1$ is conjugate to
a surface invariant by parabolic translations, see \cite{[Dan]} or
\cite{[Earp]}. Furthermore any helicoid with pitch $\ell >1$ is
conjugate to a minimal surface invariant by hyperbolic translations,
see \cite{[Earp]}.
\end{exa}

\medskip

\begin{rem}\label{R.asso}
 Assume that $\M=\m2$ and let us consider
$X,X^* : \Omega\rightarrow \m2\times \R$
two conformal minimal immersions. Let $(g,\eta)$
(resp. $(g^*,\eta^*)$) be the Weierstrass representation of
$X$ (resp. $X^*$). We know that $X$ and $X^*$ are associate in the usual
meaning, that is in the
Euclidean space $\m 3$, if and only if $g^*=g$ and
$\eta^*=e^{i\theta}\eta$ for
a real number $\theta$. Let us set $X=(h, f )$ and $X^*=( h^* , f^* )$. Since
$Q(h)=-\frac{(\eta)^2}{4}$ we see that if $X$ and $X^*$ are associate in
the usual meaning then there are associate in the meaning of definition
\ref{D.asso}. Conversely, assume that $X$ and $X^*$ are  associate in the sense
of definition \ref{D.asso}, so we have
$\eta^*=\pm \eta$ and $\vert g^*\vert=\vert g\vert$
or $\vert g^*\vert=1/\vert g\vert$.
Therefore there exists an isometry $\Gamma$ of $\m 3$ such that $X$
and $\Gamma \circ X^*$ are associate in the usual meaning.

 For example the Weierstrass representations
$(g,\eta)$ and $(e^{i\theta}g,e^{i\theta}\eta)$ for $\theta\not= 2k\pi,\ k\in
\zed$, are associate in the sense of definition  \ref{D.asso} but are
not in the usual sense.

Therefore, in $\m2\times\R$ the two notions
of associate minimal immersions are only equivalent up to an isometry
of $\m2\times\R$.
\end{rem}

\medskip

It is known that any two isometric conformal minimal immersions in
$\m3$ are associate up to an isometry. Also, it is shown in
\cite{[S-T]} that any two
isometric screw motion minimal complete immersions in
$\hi2 \times \R$ are associate.
The following example shows that this is no
longer true for any isometric immersions in $\hd \times \R$.

\medskip

\begin{exa}
\label{Ex.nonassoc}
It is given in \cite{[S-T]} an example
of complete minimal surface in
$\hd \times \R$ with intrinsic curvature constant and equals to $-1$,
$K\equiv -1$. Namely in (48) of Corollary 21
of \cite{[S-T]} we set, keeping
the notations, $H=0,\  l=m=1$, $d$ is any positive real number, $d>0$.
Therefore $U(s)=\sqrt{1+d^2}\ch (s)$, $s \in \R$. Consequently from
Theorem 19 in \cite{[S-T]} we obtain (see (33), (36) and (37))
\begin{eqnarray*}
\rho (s)&=& \argch (\sqrt{1+d^2}\ch s)\\
\lambda (\rho(s)) &=& d\int \frac{U(s)}{U^2(s)-1}ds \\
\varphi (s,\tau) &=& \tau - d\int \frac{1}{U(s)(U^2(s)-1)}ds
\end{eqnarray*}

 Now let us consider the map
$T: \m2 \rightarrow \hd \times \R$ defined for every
$(s,\tau)\in \m2$ by
$$
T(s,\tau) =\left(\ta (\rho(s)/2) \cos \varphi(s,\tau),
\ta (\rho(s)/2) \sin \varphi(s,\tau),
\lambda (\rho (s))+\varphi (s,\tau)\right),
$$

It is shown that $T$ is a regular minimal {\em embedding} with induced
metric
$$
ds_T^2= ds^2 +U^2(s)d\tau^2= ds^2 +(1+d^2)\ch^2(s)d\tau^2.
$$
A straightforward computation shows that the intrinsic curvature is given
by $K=-U^{\prime\prime}/U$. Therefore we have $K\equiv -1$. By
construction the surface $T(\m2)$ is invariant by screw motions. The
immersion $T$ is not conformal but setting $r:=\int (1/U(s))ds$ the
new coordinates $(r,\tau)$ are conformal, that is the immersion
$\wt{T} (r,\tau) :=T(s,\tau)$ is conformal. Thus the surface $T(\m2)$
is isometric to the hyperbolic plane
$(\di,\sigma^2(z)\vert \rmd z\vert^2)$, that is there exists a conformal
minimal immersion $X:\di \rightarrow \hd \times \R$ such that the
induced metric on $\di$ is the hyperbolic one and $X(\di)=T(\m2)$.
Clearly the canonical immersion $Y:\di \rightarrow \hd\times \R$
defined by $Y(z)=(z,0)$ is isometric to $X$. According to
the remark \ref{associate},
 we deduce that $X$ and $Y$ are not associate
since the third component of the Gauss map of $X$ never is equal to
$\pm 1$ as it is the case for $Y$.

Let us observe that in \cite{[H]} it can be found other examples of
complete minimal surfaces in $\hd \times \R$ with intrinsic curvature
equal to $-1$.
\end{exa}

\bigskip

\begin{rem}\label{R.asso.1}
The second author has constructed in \cite{[Earp]} new families
of complete minimal immersions in $\hi2 \times \R$ invariant by
parabolic or hyperbolic screw motions. It is shown (see Theorem 4.1)
that
any two minimal isometric parabolic screw motion immersions into
$\hi2 \times \R$ are associate. However this is no longer true for
hyperbolic screw motion immersions: there exist
isometric minimal hyperbolic  screw motion immersions into
$\hi2 \times \R$ which are not associate, see Theorem 4.2.
\end{rem}


\begin{thebibliography}{10000}

\bibitem{[A-N]} K. Akutagawa and S. Nishikawa. {\em The Gauss map and
    spacelike surfaces with prescribed mean curvature in Minkowski 3-space};
Tohoku Math. J. {\bf 42}, 67-82 (1990).

\vskip2mm
\bibitem{alhfors}L. Ahlfors. {\em On quasi-conformal mapping};
McGraw-Hill,New York.

\bibitem{[Dan]} B. Daniel. {\em Isometric immersions into
$\Sp ^n \times \R$ and $\hi n \times \R$ and applications to minimal
surfaces}. Institut de Math. de Jussieu, Pr{\'e}publications 375, Juin 2004.

\vskip2mm

\bibitem{[Du]} P. Duren. {\em Harmonic mappings in the plane}.
Cambridge Tracts in Mathematics, 156. Cambridge University Press 2004.
\vskip2mm
\bibitem{lawson}B. Lawson. {\em Lectures on minimal submanifolds}; T1,
 Mathematics Lecture Series; 009, Publish or Perish.

\vskip2mm
\bibitem{hil} U. Dierkes, S. Hildebrandt, A. K\"uster and O. Wohlrab.{\em Minimal surfaces I and II}; A series of comprehensive studies in mathematics. Springer Verlag.

\vskip2mm


\bibitem{[H]} L. Hauswirth. {\em Generalized Riemann examples in
    three-dimensional manifolds}; to appear in Pacific Journal of
Mathematics.

\vskip2mm

\bibitem{[Ken]} K. Kenmotsu. {\em Weierstrass formula for surfaces of prescribed
 mean curvature}; Math. Ann. {\bf 245}, 89-99, 1979.

\vskip2mm

\bibitem{[M-R]} W. Meeks, III and H. Rosenberg.
{\em The theory of minimal surfaces in $M\times R$}. To appear in
Comment. Math. Helv.

\vskip2mm



\bibitem{[N-R]} B. Nelli and H. Rosenberg.
{\em Minimal surfaces in $H^2\times R$}.
Bull. Braz. Math. Soc. \textbf{33}, 263-292, 2002.

\vskip2mm


\bibitem{[R]} H. Rosenberg. {\em Minimal surfaces in $M^2\times R$}.
Illinois J. Math. \textbf{46}, 1177-1195, 2002.



\vskip2mm

\bibitem{[Earp]} R. Sa Earp. {\em Parabolic and hyperbolic screw
    motion surfaces in $\hi2 \times \R$}. Preprint.
\vskip2mm.
\bibitem{sy}R. Schoen and S. T. Yau. {\em Lectures on harmonic maps};
 Conference Proceedings and lecture Notes in Geometry and Topology, II.
International Press, Cambridge, MA, 1997.

\vskip2mm

\bibitem{[S-T]} R. Sa Earp and Eric Toubiana. {\em Screw motion
    surfaces in $\hd\times \R$ and $\sd\times \R$};
Illinois Journal of Mathematics, 2005.
\vskip2mm










\end{thebibliography}
\end{document}